\documentclass{amsart}
\usepackage{amsmath,amsfonts}
\usepackage{amscd,amsthm}
\usepackage{graphicx}

\newcommand{\R}{\mathop{\mathbb{R}}}

\newcommand{\N}{\mathop{\mathbb{N}}}
\newcommand{\C}{\mathop{\mathbb{C}}}
\newcommand{\D}{\mathop{\mathbb{D}}}

\newtheorem{theorem}{Theorem}[section]

\newtheorem{remark}[theorem]{Remark}

\numberwithin{equation}{section}

\begin{document}

\title{Normalization of Poincar\'{e} Singularities {\it via} Variation
  of Constants} 
\author{T. Carletti , A. Margheri and M. Villarini}

\date{26 July 2004}

\address[Timoteo Carletti]{Scuola Normale Superiore piazza dei Cavalieri 7,
 56126  Pisa, Italy}
\address[Alessandro Margheri]{Fac. Ci\^encias de Lisboa and Centro de 
Matem\'atica e Aplica\c c\~oes Fundamentais, Av. Prof. Gama Pinto 2, 1649-003
 Lisboa, Portugal}
\address[Massimo Villarini]{Universit\`a di Modena e Reggio Emilia, via Campi
 213/b 41100 Modena, Italy}

\email[Timoteo Carletti]{t.carletti@sns.it}
\email[Alessandro Margheri]{margheri@ptmat.fc.ul.pt}
\email[Massimo Villarini]{villarini@unimo.it}



\begin{abstract}
We present a  geometric proof of the Poincar\'e-Dulac Normalization Theorem 
for analytic vector fields with singularities of Poincar\'e type.
Our approach allows us to relate the size of the convergence domain of the 
linearizing transformation to the geometry of the complex foliation associated 
to the vector field. 

A similar construction is considered in the case of linearization of maps in 
a neighborhood of a hyperbolic fixed point.
\end{abstract}

\maketitle

\section{Introduction}

Let $X : U \subset {\C}^n \mapsto T {\C}^n$ be a vector field,
holomorphic in the domain $U$; let $o \in U$ and $X(o)=0$: $o$ is a {\it
  Poincar\'{e} singular point} of $X$ if the differential $d_{o} X$
has eigenvalues $\lambda_1, \ldots , \lambda_n$ satisfying: 
\begin{equation*}
  0 \notin \mbox{convex hull of } \lambda_1 ,\ldots ,\lambda_n \, .
\end{equation*}

This is a geometric property of the complex foliation defined by $X$,
namely any vector field $Y=gX$, $g$ germ of unity at $o$, has at $o$
a Poincar\'{e} singular point if $X$ does. The geometric content of
this condition is captured by the following remark by Arnold~\cite{ar69}
 ({\it Arnold's Transversality Condition}), which is
crucial for our normalization method. 

Let $S_R$ be an Euclidean sphere
in ${\C}^n$ of radius $R$; we say that $X$ is {\it transversal} to
$S_R$ at $p \in S_R$ if $<X(p)>^{\R} \oplus T_p S_R = (T_p {\C}^n)^{\R}$. 

\begin{theorem}[Arnold~\cite{ar69}]
\label{arnold}
Let $o$ be a Poincar\'{e} singular point of $X$: then there exists
$R_0 >0$ such that for every $0<R<R_0$, $X$ is transversal to $S_R$. 
\end{theorem}

We remark that $R_0$ depends on $X$ only through its non--linear terms, being 
a linear vector field with Poincar\'e singularity transversal to $S_R$ for all
$R>0$.

For a given choice of coordinates $z=(z_1 , \ldots ,z_n)$, $z(o)=0$, which
we shall assume from now on, we can write:
\begin{equation*}
  X(z)=Az \, \partial_z + \cdots \, ,
\end{equation*}
where $A$ is a $n \times n$ complex matrix and dots stand for
nonlinear terms. We assume that in $z$
coordinates: $A=S+\varepsilon N$, $\varepsilon > 0$, is the Jordan
decomposition 
of $A$, and moreover $S=diag (\lambda_1, \ldots , \lambda_n)$. We also
introduce, 
for later use, the vector $\underline{\lambda}=(\lambda_1, \ldots ,
\lambda_n)$. 

Given an analytic diffeomorphisms $f$ one can consider
the {\em push forward} of the vector field $X$ under $f$: 
$X_*(z)=df \cdot X(f^{-1}(z))$. Geometrically this represent the
vector field $X$ in the new coordinates system determined by $f$ and
we will say 
that $X$ and $X_*$ are analytically {\em conjugated}.
For a given $X$ a natural question is to determine the ''simplest form'' it  can assume up to analytic conjugation, or given 
a vector field $X_0$ one can be interested in determining all the vector fields
that are conjugated to it.

The most interesting case occurs when  such a simplest form 
is the linear part of the vector field at the singular point.
 It is the {\em linearization} problem: it has been considered by Poincar\'{e}
  in his thesis~\cite{poincare}, and solved by him in the
case of Poincar\'{e} singularities. His results were later
generalized by Dulac~\cite{dulac} to the {\em normalization problem}. Let us
 briefly
recall what normalizing a vector field means. By a holomorphic change
of coordinates: 
\begin{equation*}
  w = z + h(z) \, ,
\end{equation*}
we try to reduce $X(z)$ to a simplest form, possibly to the linear vector field: 
\begin{equation*}
  X_{{\it lin}}=A w \, \partial_w \, .
\end{equation*}
Obstructions to realize this program are the {\it resonances}: there exist
$\underline{m} = (m_1 , \ldots ,m_n )$, $m_l \in \N$, 
$\lvert \underline{m} \rvert = m_1 + \cdots m_n \geq 2$, and 
$j\in \{1, \dots ,n\}$, such that:
\begin{equation*}
<\underline {\lambda} , \underline {m}> - \lambda_j =0\, .
\end{equation*}
A {\em formal
change of coordinates} leads to the following {\em formal normal form(s)} for
the differential equation associated to the vector field, for all
$j\in \{1, \ldots ,n\}$: 
\begin{equation*}
  \dot w_j = (A w)_j + \sum_{\substack{\underline{m}\in \N^n \, :\, \lvert
 \underline{m}\rvert \geq 2 \\ <\underline{\lambda} ,
 \underline{m}> 
 - \lambda_j = 0}} c_{\underline {m} , j} w^{\underline{m}} \, ,
\end{equation*}
where we used the standard notation $w^{\underline {m}}=w_{1}^{m_1} \cdots
  w_{n}^{m_n}$ and $(c_{\underline {m} , j})_{\underline {m} , j}\subset \C$.

 In the case of a Poincar\'{e} singular point
  there are at most finitely many resonant terms, and the
  non--resonant terms 
are bounded from below by some universal positive constant: 
$\lvert <\underline{\lambda} , \underline{m}> - \lambda_j\rvert > c$, for all 
$\underline{m}\in \N^n$ and $j \in \{ 1,\dots,n\}$, s.t. 
$<\underline{\lambda} , \underline{m}> - \lambda_j \neq 0$. This remark prevent
 the formal normalizing method from the {\it small divisor problem};
 moreover, any normal form is in this case polynomial.

\begin{theorem}[Poincar\'{e}, Dulac]
\label{pd}
Let $X$ be holomorphic in $U$ and let $o \in U$ be a Poincar\'{e}
singular point. Let $X_0$ be a polynomial normal form of $X$. Then
there exists a neighborhood $V \subset U$ of $o$ and an holomorphic
diffeomorphisms $H$ defined in $V$ such that: 
\begin{equation*}
H_* X = X_0 \, .
\end{equation*}
\end{theorem}
Even if elementary, Poincar\'{e}'s original proof of this result, as
any other more recent proof (see {\it e.g } \cite{chow}), is not 
explicit in determining the transformation
$H$ and its convergence domain. We shall give a geometric proof of
this classical theorem, {\it via} a variation of constants
approach. This will allow us to get a more explicit definition of the
normalizing transformation, and will lead us to relate the size of the
domain of the linearizing transformation to the transversality radius
$R_0$ entering in Arnold's Transversality Condition. The key idea to obtain
this result is classical: we use Hurwitz's Theorem to prove the existence of a local
biholomorphisms; then, applying Cauchy's estimates, we can extend the domain of
injectivity to the whole domain of definition.

The method used in the proof of our main result
 is an extension of the smooth normalization argument used by 
Sternberg~\cite{sternberg}. Hence it is different from all perturbative--like, or
KAM--like methods, 
where one try to push the non linear term of the vector field, to higher
and higher order through an iterative algorithm (see Remark~\ref{rem:push}).

In section~\ref{sec:normmaps} we will consider the case of discrete
time dynamical systems, i.e. the {\em Siegel Center Problem}: linearization of
a biholomorphic map in a neighborhood of a fixed point. We will present a geometric
construction, similar to the one given for flows, which allows us to solve
the Siegel Problem in the case of Poincar\'e fixed point, i.e. hyperbolic case, 
obtaining moreover an explicit bound on the size of the convergence domain of 
the linearizing map, related to geometric properties of the orbit space of the
 biholomorphism. We will also compare our result with other classical 
ones~\cite{poincare,koenigs}.

To conclude this introduction, let us briefly mention the normalization problem 
for Siegel's singularities of analytic vector fields. In this case normalization
 is not always possible: one need some additional hypotheses on the growth rate of
the small divisors and on the geometry of the foliation associated to the resonant
normal form~\cite{bruno}.

It would be interesting to deal with the Siegel case using ideas similar to those 
introduced here, to get Bruno's results. We could not succeed in developing this 
approach due to the fundamental role played  
in our geometric normalization of Poincar\'e singularities by Arnold's Transversality Condition, which is no longer
true in the Siegel case. 

\section{Normalization {\it via} Variation of Constants}

We start the description of our approach to normalization by a slight
and straightforward generalization of the variation of constants
formula. Let $X$ be, in given $z$--coordinates, of the form: 
\begin{equation*}
X(z)= X_0 (z) + X_1 (z)\, ,
\end{equation*}
where $X$, $X_0$, $X_1$ are holomorphic in the common domain $U \in
{\C}^n$. We denote by $\Phi_{X} ^{T} (z)$ (respectively $\Phi_{X_0} ^{T} (z)$) the
complex flows of $X$ (respectively $X_0$): they are both defined in a common domain
in ${\C} \times {\C}^n$. We look for a $T$-depending holomorphic
diffeomorphism $L_{T} (z)$ such that:
\begin{equation}
\label{eq:flow}
\Phi_{X} ^{T} (z) = \Phi_{X_0} ^{T} (L_{T}(z))\, ,
\end{equation}
hence, for sufficiently small positive $\Delta$ and $R$:
\begin{equation*}
L_{T} (z) = \Phi_{X_0} ^{-T} (z) \circ \Phi_{X} ^{T} (z)
\quad \forall \, \lvert T \rvert < \Delta \, , \Vert z \Vert < R \, ,
\end{equation*}
being $\Vert z \Vert$ the Euclidean norm in ${\C}^n$. 
This is the {\it Variation of Constant Transformation}. 

Let us introduce an {\em integral representation} of
this transformation, which turns to be well--suited for our
use. Differentiating relation~\eqref{eq:flow} w.r.t. time we get:
\begin{equation*}
X_0 ( \Phi_{X} ^{T} (z) ) + X_1 ( \Phi_{X} ^{T} (z) ) = X_0 (
\Phi_{X} ^{T} (z) ) + d _{\Phi_{X} ^{T} (z)} \Phi_{X_0} ^{T}  \dot
L_{T} (z) \, ,
\end{equation*}
hence:
\begin{equation*}
d _{\Phi_{X} ^{T} (z)} \Phi_{X_0} ^{T}  \dot L_{T} (z)=X_1 (\Phi_{X}
^{T} (z))\, . 
\end{equation*}
Therefore, integrating along any smooth path joining $0$ and $T$ lying
inside the disk of radius $\Delta$ in ${\C}$, we obtain:
\begin{equation}
\label{eq:intform}
L_{T} (z) = z + \int_{0}^{T} d _{\Phi_{X} ^{s} (z)} \Phi_{X_0}
^{-s} \, X_1 (\Phi_{X} ^{s} (z)) ds \, , 
\end{equation}
which is the {\it Variation of Constants Formula}. Let us explicitly
observe that such a definition depends on $X_0$, {\it i.e.} on the
chosen normal form in our case of use. 

Another way to characterize a singular point $o$ of Poincar\'{e} type
is the following one, see also~\cite{bruno} \S II page 165. There exists a line 
$l_{\omega_0}$ in $\C$: $l_{\omega_0} =\{t {\omega}_0 + \eta : {\omega}_0 \in 
S^1 \, , \eta \in \C^* \, , t\in\R\}$, such that all the eigenvalues $\lambda_1 ,
\ldots ,\lambda_n$ of $d_o X $ are contained in the halfspace
$\mathcal{H}_{\underline{\lambda}} \subset {\C}$ defined by $l_{\omega_0}$ and 
such that $0 \notin \mathcal{H}_{\underline{\lambda}}$. Let $
{l_{\omega_0}}^{\perp}  = \{ i t {\omega}_0 : t\in \R \}$
and let $ ( {l_{\omega_0}}^{\perp} )^+$ be the ''positive'' halfline not contained
in $\mathcal{H}_{\underline{\lambda}}$.
 We remark that we can change ${\omega}_0$ into $\omega$, with
 $\lvert \arg (\omega - {\omega}_0 )\rvert  < \theta$, $\theta >0$
 sufficiently small, keeping the geometric characterization of the
Poincar\'e singularity.
For a fixed choice of $\omega_0$, let $\mu_j$ be the {\em distance} of $\lambda_j$
 from ${l_{\omega_0}}$, $j\in \{1, \ldots , n\}$,
and let $\alpha = \min \{\mu_1, \ldots , \mu_n \} $, $\beta = \max
\{\mu_1, \ldots , \mu_n \} $.

We are now able to state our main result, which is a version of Theorem~\ref{pd}
 containing a more explicit definition and a geometric estimate of the domain
 of convergence of the normalizing transformation.

\begin{theorem}\label{main}
Let $X:U \in {\C}^n \mapsto T {\C}^n $ be a holomorphic vector field and let $o$
 be a Poincar\'{e} singular point of $X$. Let $R_0$ be the transversality radius
 appearing in the Arnold's Transversality Condition: hence $B_{R_0} \subset U $. 
Let:
\begin{equation*}
X(z)=X_0 (z) + X_1 (z)\, ,
\end{equation*}
where $X_0 (z)$ is a polynomial normal form of $X$, whose linear part is
in Jordan canonical form: $S+\varepsilon N$, $\varepsilon >0$, and $X_1$ is 
$m$--flat, $m > \max \{ \deg X_0 ,  
\frac{q(\beta + \varepsilon)}{\alpha + \varepsilon } \}$, being $q>1$. Then:
\begin{equation*}
L(z) = \lim _{t \mapsto + \infty } \Big \{ z + \int_{0}^{t} 
d _{\Phi_{X}^{is \omega} (z)} \Phi_{X_0} ^{-is \omega} \, X_1 
(\Phi_{X} ^{is \omega} (z)) \, ds \Big \}\, ,
\end{equation*}
where integration is along any halfline not contained in 
$\mathcal{H}_{\underline{\lambda}}$, is a normalizing biholomorphism in
 a neighborhood of $0$:
\begin{equation*}
L_* X = X_0 \, .
\end{equation*}
Moreover in the linearizing case, i.e. when $X_0 = Az \, \partial_z$
 and it is not resonant or $X_1$ is $m$--flat, $L(z)$ 
is a linearizing biholomorphism defined in a domain containing the 
Euclidean ball $B_{R_0}$ where Arnold's Transversality Condition holds.
\end{theorem}

\begin{remark}
The statement of the above theorem applies when $X$ is in a ''prepared''
normal form $X = X_0 + X_1 $: while this condition can always be satisfied
after a $m$--degree polynomial change of coordinates, this should be
taken into account on applying the bounds on the convergence domain of
the linearization. On the other hand, such ''prepared''
normal form transformation has no influence on the estimate on $m$. 
\end{remark}

\begin{remark}
Let $X=Az \, \partial_z+X_1$ then if $X_1$ is sufficiently flat, we can 
{\em linearize} $X$ even if $A$ is resonant. We recall that the maximal modulus
of resonance in the Poincar\'e case is bounded by $\beta/\alpha$, hence giving a simple interpretation of the order $m$ of flatness appearing in the statement of the theorem.

On the other hand if the order of
$X_1$ is too
small, performing the polynomial change of variables to put $X$ into the 
''prepared normal form'', one cannot avoid the ''introduction'' of resonant 
monomials of small degree, s.t. in the ''prepared normal form'', $X_0$ will be
no longer linear and our result guarantees only normalizability.
\end{remark}

As a final remark, we observe that some explicit bound on the size of the convergence domain of the normalizing 
transformation can be obtained not only in the linearization case, but in the general case of normalization, too. This bounds, if needed {\it e. g. } in applied bifurcation problem, can be easily deduced from the following proof, but we will omit them as they have not such a synthetic and geometric interpretation as in the case of linearization.

The rest of this section is devoted to the proof of Theorem~\ref{main}: 
first we will deal with the general normalization problem, and then we will show 
how to modify the arguments in the simpler case of linearization in order to get 
the estimate on the size of the linearization domain.

Without loss of generality we suppose that $\mathcal{H}_{\underline{\lambda}}
 = \{z \in {\C} : \Re z <0\}$: all the arguments in the proof transfer 
literally to the general case just considering as integration path 
$(l_{\omega_0}^{\perp})^+$ instead of the real positive semiaxis.
From a geometric point of view this choice corresponds to a time reparametrization
of the complex foliation associated to $X$, by a complex non zero factor.
Under this assumptions we have:
\begin{equation*}
\beta = \Re \lambda_n \leq \Re \lambda_{n-1} \leq \cdots \leq \Re \lambda_1 =
 \alpha <0 \, ,
\end{equation*}
where we changed the previous definitions of $\alpha , \beta $ by switching sign:
 this has no effects on the statement of the theorem and will simplify notations.

A fundamental remark is that under these hypotheses the differential equation with
 real independent variable $t$ and complex phase space $U$ given by:
\begin{equation}
\frac{d z}{dt} = X(z)\, ,
\label{reale}
\end{equation}
defines an analytic flow $\Phi_{X}^{t}(z)$ for $z \in B_{R_0}$ and $t>0$, moreover:
\begin{equation*}
\lim_{t \mapsto + \infty} \Phi_{X}^{t}(z) = 0 \, .
\end{equation*}
Of course, this is nothing but Arnold's Transversality Condition.

We can extend this remark to obtain a kind of asymptotic stability of
the origin  
as a singular point of a differential equation defined by $X$ and with
independent 
 variable $T$ varying in a sectorial neighborhood, centered on the
 real positive  
semiaxis, of infinity in the Riemann sphere. In fact, for any $0<R^0 <
R_0 $ we can  
find $\theta >0 $ such that the equations:
\begin{equation}
\frac{d z}{dt} = i \omega X(z) \, ,
\label{reale2}
\end{equation}
where $ t \in \R $ and $\vert \arg (\omega  )\vert < \theta$
define real flows which, by the same arguments we used for the
equation~\eqref{reale}, have the origin as an asymptotically stable
stationary point. Therefore such real flows imbed into the complex
flow $\Phi_{X}^T (z)$ of $X(z)$ which turns to be defined for $(T,z)
\in {\mathcal S}  \times \{ z: \Vert z \Vert  <R^0 \}$, where
$\mathcal{S} = \{T: \vert T \vert < \Delta \} \cup \{ T= t \omega ,
\vert \arg 
(\omega  ) \vert < \theta , t > 0 \} = \D_{\Delta} \cup 
\mathcal{C}_{0,\theta}$. 

We shall prove now that, for small enough $R>0$:
\begin{equation*}
L(T,z)=:L_{T} (z) : {\mathcal S}\times \mathbb{D}_R \mapsto {\C}^n \, .
\end{equation*}
The first step is to obtain an estimate for the grow 
rate of $\Vert \Phi_{X}^t (z) \Vert$. For all $j\in \{ 1, \dots , n\}$ we have:
\begin{eqnarray}
\frac{d}{dt}\Big\lvert \left(\Phi_{X}^t (z)\right)_j \Big\rvert ^2 &=&
\left[ \frac{d}{dt}\left(\Phi_{X}^t (z)\right)_j\right] 
\overline{ \left( \Phi_{X}^t (z)\right)_j} +
\left(\Phi_{X}^t (z)\right)_j
\left[\frac{d}{dt}\overline{ \left( \Phi_{X}^t (z)\right)_j}\right]\notag \\
&=& X_j (\Phi_{X}^t (z)) \overline{ \left( \Phi_{X}^t (z)\right)_j} +
\left( \Phi_{X}^t (z)\right)_j \overline{X_j (\Phi_{X}^t (z))} \notag \\
&=& 2 \Re \lambda_j \lvert\left(\Phi_{X}^t (z)\right)_j \rvert^2 +
2 \varepsilon \Re \left( (N \Phi_{X}^t (z))_j 
\overline{ \left( \Phi_{X}^t (z)\right)_j}\right)
+ {\mathcal O}\left( \lvert\left(\Phi_{X}^t (z)\right)_j \rvert^3\right)\, .
\end{eqnarray}
Hence for all $\delta >0$ we get, for all $\Vert z \Vert <R$, $R>0$ small enough:
\begin{equation*}
\frac{d}{dt} \Vert \Phi_{X}^t (z) \Vert ^2 \leq 2 
(\alpha + \varepsilon +\delta) \Vert \Phi_{X}^t (z) \Vert ^2 \, ,
\end{equation*}
and therefore:
\begin{equation}
\Vert \Phi_{X}^t (z) \Vert ^2 \leq e^{2(\alpha + \varepsilon+\delta ) t } R^2 \, . 
\label{stima1}
\end{equation}
In order to get an estimate for the integral from of $L_{T} (z)$ given 
by~\eqref{eq:intform} we need to prove firstly that 
$d _{\Phi_{X} ^{s} (z)} \Phi_{X_0} ^{-s}$ is defined for $\Vert z \Vert < R$, 
$R$ sufficiently small, and for all $T \in {\mathcal S}$, then we must find a 
suitable asymptotic estimate of it. 

For sufficiently small $\vert T \vert$ and $\Vert z \Vert$, the couple 
$( \Phi_{X}^{T} (z), d_{\Phi_{X} ^{T} (z)} \Phi_{X_0} ^{-T})$ is the solution 
of the following Cauchy problem for a system at variation type:
\begin{equation*}
\begin{cases}
\dot w & = X(w ) \\
\dot W & = - \frac{\partial X_0}{\partial w} W \\
w (0) & = z \\
W (0) & = E \, ,
\end{cases} 
\end{equation*}
where $E$ is the identity matrix in ${\C}^n$, $w \in {\C}^n$ and 
$W \in {\C}^n \times {\C}^n$. Therefore the existence of 
$t \mapsto (w (T), W(T)) = ( \Phi_{X} ^{T} (z), d _{\Phi_{X} ^{T} (z)}
 \Phi_{X_0} ^{-T})$ for every real $t>0$ and $\Vert z \Vert < R^0 $  
follows from the asymptotic stability of the origin as a singular point 
of \eqref{reale} and from basic theory of linear ordinary differential 
equations. To get the desired asymptotic estimate for 
$d _{\Phi_{X} ^{t} (z)} \Phi_{X_0} ^{-t}$ 
we consider the above system with fixed $z$ and 
writing $X_0 (w ) = A w + g ( w )$ we obtain the equation at variation:
\begin{equation*}
\dot W = -A W - \frac{\partial g}{\partial w } (\Phi_{X} ^{t})W \, .
\end{equation*}
From this equation is readily obtained the following inequality for the norm 
of linear operators:
\begin{equation*}
\Vert d _{\Phi_{X} ^{t} (z)} \Phi_{X_0} ^{-t}) \Vert \leq \Vert e^{-t A} \Vert + 
 {\mathcal O} (\Vert \Phi_{X} ^{t} \Vert)\, ,
\end{equation*}
and then for every $0<R<R^0$ there exits $\delta >0$  such that for every $t>0$:
\begin{equation}
\label{stima2}
\Vert d _{\Phi_{X} ^{t} (z)} \Phi_{X_0} ^{-t}) \Vert \leq 
e^{-  (\beta +\varepsilon- \delta ) t } R\, .
\end{equation}
We can give now a uniform bound on $\Vert L_t (z) \Vert$ for $t>0$ and 
$\Vert z \Vert < R$. Recalling that $\Vert X_1 (w ) \Vert \leq C \Vert w \Vert ^m $ 
and $ m > \frac{q(\beta +\varepsilon )}{\alpha +\varepsilon}$, so that 
$(\beta +\varepsilon - \delta ) - m (\alpha + \varepsilon +\delta ) >0$ for 
sufficiently small $\delta$, from the estimates~\eqref{stima1} and~\eqref{stima2}
 we get:
\begin{eqnarray}
\Vert L_t (z) - z \Vert &=& \Big\Vert \int_{0}^{t} d_{\Phi_{X}^{s} (z)} 
\Phi_{X_0}^{-s} X_1 (\Phi_{X}^{s} (z)) \, ds \Big\Vert 
\leq C \Vert z \Vert ^m \int_{0}^{t} e^{[-(\beta +\varepsilon - \delta ) + 
 m (\alpha +\varepsilon+ \delta )] s} \, ds \notag \\
&\leq & \frac{C}{(\beta +\varepsilon- \delta )-m (\alpha +\varepsilon + \delta )}
 R^m (1 - e^{-[(\beta +\varepsilon - \delta )+m (\alpha + \varepsilon +\delta )]t}) \, ,
\end{eqnarray}
for $t>0$ and $\Vert z \Vert < R$. Hence every $L_t $, $t>0$, maps the Euclidean 
ball $B_R$ into $B_{R + R^{\prime}}$ where $R^{\prime} =
 \frac{C}{(\beta +\varepsilon - \delta ) - m (\alpha +\varepsilon + \delta )} R^m$.

 An analogous estimate leads to the proof of the convergence of $t \mapsto L_t $. 
In fact, let us suppose for the moment that $\tau $ is real and with sufficiently 
small modulus, then:
\begin{eqnarray}
\Vert L_{t + \tau }(z) - L_t (z) \Vert &=& \Big \Vert  \int_{t}^{t + \tau} 
d_{\Phi_{X}^{s} (z)} \Phi_{X_0}^{-s} X_1 (\Phi_{X} ^{s} (z)) \, ds \Big \Vert \notag \\
&\leq & C \Vert z \Vert ^M \frac{e^{-[(\beta +\varepsilon - \delta ) - 
m (\alpha +\varepsilon + \delta )] t}}{[(\beta +\varepsilon - \delta ) - 
m (\alpha +\varepsilon + \delta )]} (1- e^{-[(\beta +\varepsilon - \delta ) - 
m (\alpha +\varepsilon + \delta )] \tau })\, .
\end{eqnarray}
Therefore from the Cauchy condition and the hypothesis on $m$, it follows that:
\begin{equation*}
\lim_{t \mapsto + \infty} L_t (z) = L(z) \, ,
\end{equation*}
uniformly when $\Vert z \Vert < R$. By the same argument with obvious 
modifications we get:
\begin{equation*}
\lim_{\substack{T \mapsto  \infty \\T \in {\mathcal S}}} L_T (z) = L(z) \, ,
\end{equation*}
uniformly when $\Vert z \Vert < R$.

To end the proof of existence of a locally defined normalizing transformation 
we need to show that $L(z)$  conjugates the vector field to $X_0$ i.e. : $L_* X =
X_0$ and is a biholomorphism in a neighborhood of the origin.

To prove the first claim is enough to show that $L$ conjugates the corresponding 
flows, namely:
\begin{equation*}
\Phi_{X_0} ^{- \tau} \circ L \circ \Phi_{X} ^{\tau} = L \quad \forall \tau\, .
\end{equation*}
This is obvious writing:
\begin{equation}
L= \lim_{\substack{T \mapsto  \infty \\T \in {\mathcal S}}} \Phi_{X_0}^{- t}
 \circ \Phi_{X}^{t} \, ,
\label{coniugio}
\end{equation}
as in this case:
\begin{equation*}
\Phi_{X_0}^{- \tau} \circ L \circ \Phi_{X}^{\tau} = 
\lim_{\substack{T \mapsto  \infty \\T \in {\mathcal S}}} \Phi_{X_0}^{- (t + \tau ) }
  \circ \Phi_{X} ^{t + \tau } = L(z)\, .
\end{equation*}
So the first claim is proved if~\eqref{coniugio} holds; let us prove it. Because:
\begin{equation*}
L_{T} (z) = \Phi_{X_0}^{- T} \circ \Phi_{X}^{T}(z) = z + \int_{0}^{T} 
d_{\Phi_{X}^{s} (z)} \Phi_{X_0}^{-s} X_1 (\Phi_{X}^{s} (z)) \, ds \, ,
\end{equation*}
for sufficiently small $\vert T \vert $ and $\Vert z \Vert < R$, let us define 
$t_0 = \sup \{ t>0 : \mbox{ for every } \Vert z \Vert < R \mbox{ and } 
\tau \in [0, t) \smallskip : \smallskip \Vert \Phi_{X_0}^{- \tau }
\circ \Phi_{X}^{\tau }(z) 
 \Vert < + \infty  \} $, and let us suppose by contradiction that $t_0
 < \infty$. 
 Then there exists a sequence $(t_m , z^{(m)})$ such that $t_m \mapsto t_0$, 
$\Vert z^{(m)} \Vert <R$ and:
\begin{equation*}
\lim_{m \mapsto \infty } \Big\Vert  z^{(m)} + \int_{0}^{t_m} 
d_{\Phi_{X}^{s} (z)} \Phi_{X_0}^{-s} X_1 (\Phi_{X}^{s} (z)) \, ds
\Big\Vert =+\infty \, . 
\end{equation*}
This contradicts the bound $L_t (B_R) \subset B_{R+R^{\prime}}$,
 $t>0$, from which the claim follows. The proof 
 that $L$ is a biholomorphism, locally invertible in a neighborhood of
 the origin,  
follows from Weierstrass' Theorem applied to the family of analytic maps 
$\{ L_t \}$ and from the Inverse Function Theorem together with the remark that
 $d_0 L = \mbox { identity }$. 

This ends the proof of the existence of the 
normalizing transformation.

Let us come \label{page:lienar1} now to the case when $X_0 = A$ and non--resonant (or $X_1$
  is sufficiently flat), {\it i.e.} the case  when the linearizing map is:
\begin{equation}
\label{eq:linearizing}
L_{T} (z) = z + \int_{0}^{T} e^{-sA} X_1 (\Phi_{X} ^{s} (z)) \, ds \, .
\end{equation}
We must prove that $L(z) $ is a linearizing biholomorphism defined in 
$B_{R_0}$: this will follows from the proof of the analogous claim for $B_{R^0}$.

Firstly we observe that as $e^{-s A}$ is globally defined and $\Phi_{X}^{s}(z)$
 is defined for $\Vert z \Vert < R^0 $ and $ T \in {\mathcal S} $, $L_T (z)$ 
is well--defined in $B_{R^0} \times {\mathcal S}$. Moreover in the same domain 
we get the estimate:
\begin{eqnarray}
\Vert L_t (z) - z \Vert &=& \Big\Vert \int_{0}^{T} e^{-sA} 
X_1 (\Phi_{X} ^{s} (z))\, ds \Big \Vert \notag \\
& < & C^{\prime}  
+ \frac{C}{(\beta +\varepsilon - \delta ) - m (\alpha +\varepsilon+ \delta )} 
R^m\, ,
\end{eqnarray}
where $\Vert X_1 ( w ) \Vert  < C \Vert w \Vert ^m$ for $\Vert w \Vert < R $ and:
\begin{equation*}
C^{\prime} = \Big \Vert \int_{0}^{s_0}e^{-sA} X_1 (\Phi_{X} ^{s} (z)) 
\, ds \Big\Vert \, ,
\end{equation*}
for some $s_0$. Therefore each $L_T$, $(z,T) \in  B_{R^0} \times {\mathcal S}$, maps 
$B_{R^0}$ into $B_{R^0 + C^{\prime}+ R^{\prime}}$, where 
$R^{\prime}= \frac{C}{(\beta +\varepsilon - \delta ) -m (\alpha +\varepsilon 
+ \delta )} R^m$. With a similar argument we get that, if $\Re T > s_0 $:
\begin{eqnarray}
\Vert L_{T + \tau  }(z) - L_t (z) \Vert & < & \Big\Vert \int_{T}^{T + \tau }
e^{-sA} X_1 (\Phi_{X} ^{s} (z)) \,ds \Big\Vert \notag \\
& < & C R^m \frac{e^{-s_0 [(\beta +\varepsilon - \delta) - m 
(\alpha +\varepsilon + \delta )]}}{(\beta +\varepsilon - \delta) -  m 
(\alpha +\varepsilon + \delta )}\, ,
\end{eqnarray}
therefore:
\begin{equation*}
\lim_{t \mapsto + \infty} L_t (z) = L(z)\, .
\end{equation*}
Hence we have a family $\{ L_t \}_{t>0}$ of biholomorphisms from 
$B_{R^0}$ to $B_{R^0 + C^{\prime}+ R^{\prime}}$, converging in $B_{R^0}$ to $L$:
we will adapt a classical argument~\cite{cartan} concerning sequences of 
automorphisms to prove that $L$ is a biholomorphism on $B_{R^0}$, too. 
Let us denote $J_{L_t}$ and $J_L$ respectively the jacobians of $L_t$ and $L$: 
of course:
\begin{equation*}
\lim_{t \mapsto \infty} J_{L_t} = J_L \, ,
\end{equation*}
uniformly on compact subsets of $B_{R^0}$. It is a straightforward consequence of 
Hurwitz's Theorem and of the equality $J_L(0) = 1$ that $J_{L} (z) \neq 0$
 for every $z \in B_{R^0}$. Therefore the proof that $L$ is a biholomorphism from 
$B_{R^0}$ onto its image will end if we prove that $L$ is injective in $B_{R^0}$.

Let us suppose, by contradiction, that there exist 
$z^{(1)} , z^{(2)} \in B_{R^0}$ such that $L(z^{(j)})=w$, $j=1,2$.
 Let $B_r (z^{(j)})$, $j=1,2$ be two Euclidean balls centered at $z^{(j)}$ and 
having the same radius $r$, such that:
\begin{equation*}
B_r (z^{(j)}) \subset B_{R^0}\quad \text{and}\quad
B_r (z^{(1)}) \cap B_r (z^{(2)}) = \emptyset \, .
\end{equation*}
We claim that there exists $R>0$, depending on $\vert J_L (z^{(j)}) \vert$, 
$R^0$ and $R'$, and $t_0 >0$ such that for every $t>t_0 $, $j=1,2$, we have:
\begin{equation}
B_R (w) \subset L_t ( B_r (z^{(j)} ) \, .
\label{auto}
\end{equation}
This leads to a contradiction: in fact from~\eqref{auto} it follows that for 
sufficiently large $t>t_0$ there exists two points $w^{(1)}$, 
$w^{(2)}$, $w^{(j)} \in B_r (z^{(j)})$, such that 
$L_t (w^{(j)}) = w$, which is impossible because $L_t $ is a biholomorphism 
on $B_{R^0}$. Let us prove~\eqref{auto}. For $j=1,2$, there exists $t_0 > 0$ 
such that if $t>t_0$ then 
$\vert J_{L_t} (w^{(j)}) \vert \geq \frac{1}{2} 
\vert J_{L} (w^{(j)}) \vert >0$. From Cauchy inequalities we get, for $t>t_0$,
 $l,k=1, \ldots ,n$:
\begin{equation*}
\Big\vert \frac{\partial (L_t )_l }{\partial z_k } \Big\vert < \frac{R'}{r}\, ,
\end{equation*}
and therefore $\Vert d_{w} (L_t)^{-1} \Vert > \sigma \Vert w \Vert $, 
where $\sigma >0$ depends on $R',R^0 , r,n, \vert J_l (\omega^{(j)}) \vert$ but 
is independent of $t$, for $t>t_0$. Hence
\begin{equation}
\Vert d_{z^{(j)}} (L_t) \Vert < \sigma \Vert w \Vert \, ,
\label{stima}
\end{equation}
for every $w \in {\bf C}^n$. Another application of Cauchy inequalities leads 
to the following estimates of the error made substituting the linear approximation 
to the complete Taylor series, holding true when 
$\Vert z - z^{(j)} \Vert < \frac{r}{n}$:
\begin{eqnarray}
\Big \vert \sum_{ k \in {\bf N}^n : \vert  k \vert \geq 2 }
 \frac{1}{k !} D^{k } (L_t )_l (z-z^{(j)} )^{k} \Big\vert 
&\leq & R'\Big \{ \frac{1}{r^2} n^2 \Vert z - z^{(j)} \Vert ^2 + 
\frac{1}{r^3} \Vert z - z^{(j)} \Vert ^3 + \cdots \Big\} \notag \\
&\leq &\frac{r R'n^2 \Vert z - z^{(j)} \Vert ^2 }{r^2 (r-n \Vert z - z^{(j)} \Vert )} \, .
\end{eqnarray}
Therefore there exists $\delta >0$, depending on $n,R',r,\vert J_L (z^{(j)}) \vert$,
 but independent of $t>t_0 $, such that if $\Vert z - z^{(j)} \Vert < \delta $ 
then:
\begin{equation*}
 \Vert L_t (z) - L_t (z^{(j)}) - d_{z^{(j)}} L_t (z- z^{(j)}) \Vert 
< \frac{1}{2} \sigma \Vert z - z^{(j)} \Vert \, .
\end{equation*}
Combining this inequality and \eqref{stima} we obtain:
\begin{equation*}
\Vert L_t (z) - L_t (z^{(j)} \Vert > \frac{1}{2} \sigma \Vert z - z^{(j)} \Vert \, ,
\end{equation*}
hence for every $t>t_0$, $j=1,2$:
\begin{equation*}
B_{\frac{1}{3} \sigma \delta }(z^{(j)} )  \subset L_t ( B_r (z^{(j)} )\, ,
\end{equation*}
and~\eqref{auto} follows as a consequence of the convergence of the $L_t$'s to $L$. 
The proof is concluded by remarking that the conjugacy functional equation 
locally satisfied by $L$ extends to  $B_{R^0}$ by analytic \label{page:lienar2} 
continuation.

\begin{remark}
  \label{rem:push}
Let us show that for any finite $t$, the transformed vector field $(L_t)_*X$ 
has non--linearities of the same order of $X$. This is completely different from 
classical methods of perturbative 
theory where one looks for diffeomorphisms $\phi$, such that $\phi_*X$ has non 
linear terms of order higher than $X$.

To simplify assume $A$ to be non resonant and let us write 
$X_1=\sum_{\vert\underline{k}\vert = m}X_{1,\underline{k}}z^{\underline{k}}
 \partial_z +
\mathcal{O}(\vert z \vert^p)$, $p>m$. Fix $t>0$ and look for 
$h_t(z)$, s.t.: $L_t(z)=z+h_t(z)+\dots$. From~\eqref{eq:linearizing} we get:
\begin{equation*}
(h_t(z))_j=\sum_{\vert\underline{k}\vert = m}
\frac{e^{t\left(<\underline{\lambda},\underline{k} > -\lambda_j \right)}-1}
{<\underline{\lambda},\underline{k} > - \lambda_j}X_{1,\underline{k},j}
z^{\underline{k}} \, ,
\end{equation*}
where $X_{1,\underline{k},j}$ is the $j$--th component of $X_{1,\underline{k}}$.
Hence we get:
\begin{equation*}
(L_t)_*(Az\partial_z+X_1(z))= Az\partial_z
+\sum_{j\in \{1,\dots,n\},\vert\underline{k}\vert = m}
e^{t\left(<\underline{\lambda},\underline{k} > -\lambda_j \right)}
X_{1,\underline{k},j}z^{\underline{k}} \partial_{z_j}+\dots \, .
\end{equation*}
\end{remark}

\section{Linearization of biholomorphic Maps}
\label{sec:normmaps}

In this section we study the case of discrete time systems:
iteration of biholomorphisms, with particular interest in the problem of 
their linearization. We will present a modified version of the construction
 previously given for flows, which allows us to solve the linearization problem
in the case of hyperbolic fixed point.

Let us consider an analytic diffeomorphisms of $n$ complex
variables fixing the origin: $F\in {\it Diff}^{\omega}(\C^n,0)$,
$F(0)=0$ and assume that in 
the chosen coordinates $z$, it has the form: $F(z)=Az+F_1(z)$, where
$A=S+\varepsilon N$, $\varepsilon > 0$, is the Jordan Canonical form of
$dF_0$, whereas 
$\Vert F_1(z)\Vert = \mathcal{O}(\Vert z \Vert^2)$.

 Let $\mu_1,\dots,\mu_n$ be its eigenvalues, we will assume that
the origin is a fixed point of 
{\em Poincar\'e type}~\footnote{If $\inf_j \lvert \mu_j \rvert > 1$, we will
consider the map $\tilde{F}=F^{-1}$.}, namely:
\begin{equation}
  \label{eq:poincmap}
  \sup_j \lvert \mu_j \rvert < 1 \, ,
\end{equation}
and let us also introduce the vector
$\underline{\mu}=(\mu_1,\dots,\mu_n)$. 

Once again we are interested in the possibility of ''reducing'' the
given system to a ''simplest form'', $F_0(z)$, through an analytical
change of variables $H(z)$ which locally conjugates $F$ and $F_0$:
\begin{equation}
  \label{eq:conjugmap}
  H\circ F = F_0 \circ H \, .
\end{equation}

The ''simplest form'' could be the linear map $F_{\it lin}(z) =Az$,
but {\em resonances}~\footnote{In the discrete time case resonances
  are couples $j \in \{1,\dots,n\}$, $\underline{m}\in \N^n$, s.t. $\vert
\underline{m}\vert \geq 2$ and
$\underline{\mu}^{\underline{m}}-\mu_j=0$.} can be an
obstruction. In presence of resonances the ''simplest form'' is given
by the following {\em normal form}, for all $j\in \{ 1,\dots,n \}$:
\begin{equation}
  \label{eq:normformap}
  (F_0(z))_j=(Az)_j+\sum_{\substack{\underline{m}\in \N^n \, :\, \lvert
 \underline{m}\rvert \geq 2 \\ \underline{\mu}^{\underline{m}}-
 \mu_j = 0}} b_{\underline {m} , j} w^{\underline{m}} \, ,
\end{equation}
where $(b_{\underline {m} , j})_{\underline {m} , j}\subset \C$.

The main result of this section is

\begin{theorem}\label{main2}
Let $F:U \subset \C^n \mapsto \C^n $ be a biholomorphisms fixing the origin, 
assume moreover the origin to be a fixed point of non resonant Poincar\'{e} type .
 Let the chosen coordinates such that $dF_0=S+\varepsilon N$, $\varepsilon >0$, is
in Jordan canonical form, hence $F(z)=(S+\varepsilon N)z+F_1(z)$ and $F_1$ is 
$m$--flat, $m > q\frac{\lvert \log \min |\mu_j|\rvert}
{\lvert \log (\max |\mu_j|+\varepsilon)\rvert}$, $q>1$. Then:

\begin{equation*}
L(z) = z+\sum_{l=0}^{+\infty}\Delta_l(z)\, ,
\end{equation*}
where $\Delta_l(z)=(S+\varepsilon N)^{-l}f_1\circ F^l(z)$ and 
$f_1 = (S+\varepsilon N)^{-1}F_1$, is a linearizing biholomorphism in a 
neighborhood of $0$:
\begin{equation*}
AL = L\circ F\, .
\end{equation*}
The linearizing map has a convergence domain containing an euclidean ball 
of radius $R_{\delta}$ explicitly estimated in the proof by \eqref{eq:contract}.
\end{theorem}

The rest of the section is devoted to prove this result. Let us assume $A$ to 
be non--resonant, hence the normal
form reduces to the linear part of $F$.
We are looking for a family of maps $L_l(z)$, defined for $l\in \N$ and $\Vert z
\Vert$  sufficiently small, such that: $F^{\circ l}(z)=A^lL_l(z)$, 
$A=(S+\varepsilon N)$. Let us introduce the ''one--step'' map:
$\Delta_l(z)=L_{l+1}(z)-L_l(z)$. Under our assumptions we obtain:
\begin{equation}
  \label{eq:onestepmap}
  \Delta_l(z)=A^{-l}f_1\circ F^{\circ l}(z) \, ,
\end{equation}
where $f_1 = A^{-1}F_1$.

For all positive integer $l$, we trivially have the following properties:
\begin{enumerate}
\item $L_l(0)=0$;
\item $(dL_l)_0=identity$;
\item $L_{l+1}=z+\sum_{k=0}^l\Delta_k(z)$.
\end{enumerate}
From the existence of $L(z)=\lim_{l\rightarrow
  +\infty}L_l(z)$, we also get:
\begin{equation*}
  AL = L\circ F \, ,
\end{equation*}
namely $L$ linearizes $F$.

Let us now prove the existence of the previous limit.
Because there are no resonances we can perform a polynomial change of
coordinates such that $F$ is in some ''prepared form'' where $F_1(z)$
has order $m$, with $m$ arbitrary large. We remark that $\Delta_l(z)$
and $F_1(z)$ have the same order.

Let us call $\rho^*=\max_j |\mu_j|$ and $\rho_*=\min_j
|\mu_j|$, then for any $\delta >0$ we 
can find $R_{\delta}>0$ such that:
\begin{equation}
  \label{eq:contract}
  \Vert F(z)\Vert \leq (\rho^* + \varepsilon +\delta)
 \Vert z \Vert \quad \forall \Vert z
  \Vert < R_{\delta} \, .
\end{equation}
This is in some sense the analogous of the Arnold Transversality
Condition: it ensures that orbits of $F$ intersect transversally
 (in fact enter into) the euclidean ball of radius $R_{\delta}$.

By assumption $\rho^*<1$ (Poincar\'e case), hence we can
choose $\delta >0$ such that: $\rho^* + \varepsilon + \delta <1$. Using
the $m$--flatness of $f_1$ we claim that there exists a positive constant $C$ s.t.:
\begin{equation}
  \label{eq:estimatemap}
 \Vert A^{-l}f_1 \circ F^{\circ l}(z) \Vert \leq C \rho_*^{-l}
 (\rho^* + \varepsilon + \delta)^{lm} \Vert z \Vert^m \, .
\end{equation}
Because $m\geq q\lvert\log \rho_* \rvert/\lvert \log (\rho^*+\varepsilon)\rvert$, we
have: $(\rho^* + \varepsilon + \delta)^m/\rho_* = \vartheta<1$ for $\delta$ small
enough. Then:
\begin{equation}
  \label{eq:estimate2}
  \Vert \Delta_l(z)\Vert \leq C \vartheta^l \Vert z \Vert^m \, ,
\end{equation}
for all $\Vert z \Vert < R_{\delta}$.

The existence of the limit for $L_l(z)$ follows by Cauchy criterium and the 
estimate:
\begin{equation}
  \label{eq:estiamte3}
  \Vert L_{l+k}(z)-L_l(z)\Vert \leq \Big \Vert\sum_{p=l+1}^{l+k}\Delta_p(z)\Big\Vert
  \leq C \Vert z \Vert^m \frac{\vartheta^{l+1}(1-\vartheta^k)}{1-\vartheta} \, .
\end{equation}
The proof that $L$ is a biholomorphism from $R_{\delta}$ onto its image follows 
the same lines as the analogous result in the case of vector fields. And this
concludes the proof.

\begin{remark}
This linearization procedure is new and it is different from the classical ones 
of Poincar\'e~\cite{poincare} or Koenigs~\cite{koenigs}, in fact for any finite $l$,
$L_l(z)$ doesn't ''push'' the non linearities of the given biholomorphisms to higher
and higher orders. This construction is also different from the C\'esaro mean, thanks
to the presence of the term $f_1$, which also increases the speed of the
convergence.
\end{remark}

\end{document}